# Nonparametric statistics on manifolds with applications to shape spaces[*]

Abhishek Bhattacharya[1] and Rabi Bhattacharya[1]

*University of Arizona*

**Abstract:** This article presents certain recent methodologies and some new results for the statistical analysis of probability distributions on manifolds. An important example considered in some detail here is the 2-D shape space of k-ads, comprising all configurations of $k$ planar landmarks ($k > 2$)-modulo translation, scaling and rotation.

## Contents

## 1. Introduction

The statistical analysis of shape distributions based on random samples is important in many areas such as morphometrics (discrimination and classification of biological shapes), medical diagnostics (detection of change or deformation of shapes in some organs due to some disease, for example) and machine vision (e.g., digital recording and analysis based on planar views of 3-D objects). Among the pioneers on foundational studies leading to such applications, we mention Kendall [20] (also see Kendall et al. [21]) and Bookstein [9]. The geometries of the spaces are those of differentiable manifolds often with appropriate Riemannian structures.

---

[*]Supported by NSF Grant DMS 04-06-143.
[1]Department of Mathematics, University of Arizona, Tucson, AZ 85721, USA, e-mail: abhishek@math.arizona.edu; rabi@math.arizona.edu

Our goal in this article is to establish some general principles for nonparametric statistical analysis on such manifolds and apply those to some shape spaces, especially Kendall's two-dimensional shape space $\Sigma_2^k$ of the so-called k-ads, i.e., the space of configurations of $k$ points on the plane (not all identical), identified modulo size and under Euclidean motions of translation and rotation. Two sample tests for the comparison of both extrinsic and intrinsic Fréchet mean shapes and mean variations of two distributions on $\Sigma_2^k$ are provided. As far as we know the explicit computations of these tests are new. In the case of the intrinsic mean and variation, the usual support criterion (see, e.g., Le [24] and Bhattacharya and Patrangenaru [6, 7, 8]) is significantly relaxed, thereby substantially enhancing the applicability of the tests.

For recent results on statistical analysis of 3-D shapes, which we do not consider here, we refer to Dryden et al.[11] and Bandulasiri et al. [2].

Sometimes the sample sizes in shape analysis are only moderately large. Under such circumstances, one may more effectively use Effron's bootstrap methods (Effron [14]), whose superiority over the classical CLT-based confidence regions and tests may be established via higher order asymptotics (see, e.g., Babu and Singh [1], Bhattacharya and Qumsiyeh [5], Bhattacharya and Ghosh [4], Ghosh [16], Hall [17]).

We next turn to the specific example of main interest to us, namely, $\Sigma_2^k$. For purposes of medical diagnostics, classification of biological species, etc., one may use expert help to choose a suitable ordered set of $k$ points or landmarks in the plane, or a *k-ad*,

$$\mathbf{z} = \{(x_j, y_j), 1 \leq j \leq k\},$$

on a two-dimensional image of an object under consideration. One assumes that not all $k$ points are the same, and $k > 2$. Kendall's shape space $\Sigma_2^k$ comprises the equivalence classes of all such k-ads under translation, rotation and scaling. For a given k-ad $\mathbf{z}$, the effect of translation is removed by considering $\mathbf{z} - \langle \mathbf{z} \rangle$ where $\langle \mathbf{z} \rangle$ is the vector whose elements are all equal to the mean location of the k-ad, namely, $(1/k) \sum_{j=1}^{k} (x_j, y_j)$. The translated k-ads then lie in the $(2k-2)$-dimensional hyperplane $H$ of $(\Re^2)^k \approx \Re^{2k}$, given by

$$H = \{(x_j, y_j)_{1 \leq j \leq k} : \sum x_j = 0, \ \sum y_j = 0\},$$

and they comprise all of $H$ except the origin. The effect of scale, or length, is removed by dividing $\mathbf{z} - \langle \mathbf{z} \rangle$ by $\|\mathbf{z} - \langle \mathbf{z} \rangle\|$ where $\|.\|$ is the usual Euclidean norm in $(\Re^2)^k$,

$$\|(u_j, v_j)_{1 \leq j \leq k}\| = [\sum (u_j^2 + v_j^2)]^{1/2}.$$

The resulting transformed k-ad $\mathbf{w} = (\mathbf{z} - \langle \mathbf{z} \rangle) / \|\mathbf{z} - \langle \mathbf{z} \rangle\|$ is called the *preshape* of the k-ad $\mathbf{z}$. The set of preshapes is then naturally identified with the unit sphere in $H$, which is basically the same as the unit sphere $S^{2k-3}$ in $\Re^{2k-2}$. Finally, the *shape* $[\mathbf{z}]$ of a k-ad $\mathbf{z}$ is given by the *orbit* of $\mathbf{w} = (u_j, v_j)'_{1 \leq j \leq k}$ under rotation, namely,

$$(1.1) \qquad [\mathbf{z}] = \left[ \begin{pmatrix} \cos\theta & -\sin\theta \\ \sin\theta & \cos\theta \end{pmatrix} \begin{pmatrix} u_j \\ v_j \end{pmatrix}_{1 \leq j \leq k} \right], \ -\pi < \theta \leq \pi.$$

Thus $\Sigma_2^k$ is a *quotient space* of $S^{2k-3}$, namely, $S^{2k-3}/S^1$, and it has dimension $2k - 4$.

We will use a mathematically more convenient way of describing $\Sigma_2^k$ as achieved by viewing a k-ad as an element of $\mathbb{C}^k$, namely, $\mathbf{z} = (x_j + iy_j)_{1 \leq j \leq k}$. Then $\langle \mathbf{z} \rangle$ is the complex $k$-vector whose elements are all equal to $(1/k)\sum_{j=1}^{k}(x_j + iy_j)$. The translated k-ad then lies in the complex $(k-1)$-dimensional hyperplane

$$H^{k-1} = \{(a_j)_{1 \leq j \leq k} \in \mathbb{C}^k : \sum_{j=1}^{k} a_j = 0\}.$$

The norm $\|\mathbf{z} - \langle \mathbf{z}\rangle\|$ has the same value as before. But the rotation by an angle $\theta$ of $\mathbf{w} = (\mathbf{z} - \langle\mathbf{z}\rangle)/\|\mathbf{z} - \langle\mathbf{z}\rangle\|$ may now be expressed as $e^{i\theta}\mathbf{w}$. For a system of coordinate neighborhoods, or parametrization, of this *spherical representation* of $\Sigma_2^k$ as a quotient space of $S^{2k-3}$, see Gallot et al. ([15], pp. 32, 34).

Another parametrization of $\Sigma_2^k$, compatible with the above, is obtained by viewing the shape of a k-ad $\mathbf{z} \equiv (x_j + iy_j)_{1 \leq j \leq k}$ as the orbit

$$\{z_0(\mathbf{z} - \langle\mathbf{z}\rangle) : z_0 \in \mathbb{C} \setminus \{0\}\}.$$

Note that $z_0 = \lambda e^{i\theta}$ for $\lambda = |z_0|$ and some $\theta \in (-\pi, \pi]$, so that the orbit, namely, a complex line through the origin in $H^{k-1}$, is independent of both scale and rotation and, therefore, a representation of the shape of $\mathbf{z}$. Thus $\Sigma_2^k$ is (isomorphic to) the space of all complex lines through the origin in $\mathbb{C}^{k-1}$, the *complex projective space* $\mathbb{C}P^{k-2}$, a familiar and important example in differential geometry. For a system of coordinate neighborhoods for $\Sigma_2^k$ viewed as $\mathbb{C}P^{k-2}$, see Gallot et al. ([15], pp. 9, 10, 64, 65).

We next consider an *extrinsic distance* on $\Sigma_2^k$ corresponding to a special embedding, namely, the *Veronese-Whitney embedding* $\phi_E$ of $\Sigma_2^k$ into the space $S(k, \mathbb{C})$ of $k \times k$ complex Hermitian matrices:

(1.2) $$\phi_E([\mathbf{z}]) = \mathbf{ww}^*$$

where $\mathbf{w} = (\mathbf{z} - \langle\mathbf{z}\rangle)/\|\mathbf{z} - \langle\mathbf{z}\rangle\|$ is the preshape of $\mathbf{z}$. Here $\mathbf{w}$ is regarded as a column vector of $k$ complex numbers, $\mathbf{w} = (w_1, w_2, \ldots, w_k)'$, and $\mathbf{w}^*$ is the transpose of its complex conjugate. Observe that the right side of (1.2) is constant on the orbit $\{e^{i\theta}w : -\pi < \theta \leq \pi\}$ of the preshape $w$ and is, therefore, a function of the shape $[z]$ of the k-ad. Also, this function is one-to-one on $\Sigma_2^k$ into $S(k, \mathbb{C})$. The vector space $S(k, \mathbb{C})$, with the real scaler field $\Re$, has dimension $k^2$. This is because a $k \times k$ Hermitian matrix is specified by $k$ real numbers on the diagonal and $\binom{k}{2}$ complex numbers (i.e., $2\binom{k}{2}$ real numbers) as lower-right off-diagonal elements. On $S(k, \mathbb{C})$ define the norm $\|.\|$ and distance $d$ by

$$\|A\|^2 = \text{Trace} AA^* = \text{Trace} A^2,$$
(1.3) $$d^2(A, B) = \|A - B\|^2 = \text{Trace}(A - B)^2.$$

Note that this is the same as the Euclidean norm and distance in $\Re^{2k^2}$. The induced distance $\rho_E$ on $\Sigma_2^k$ is then given by

$$\rho_E^2([\mathbf{z}], [\mathbf{w}]) = d^2(\phi_E([\mathbf{z}]), \phi_E([\mathbf{w}])) = \text{Trace}(\mathbf{uu}^* - \mathbf{vv}^*)^2$$
$$= \sum_{j=1}^{k}|u_j|^2 + \sum_{j=1}^{k}|v_j|^2 - \sum_{j=1}^{k}\sum_{j'=1}^{k}(u_j \bar{u}_{j'} v_{j'} \bar{v}_j + v_j \bar{v}_{j'} u_{j'} \bar{u}_j)$$
(1.4) $$= 2 - 2|\mathbf{u}^*\mathbf{v}|^2$$

where $\mathbf{u}$ and $\mathbf{v}$ are the preshapes of $[\mathbf{z}]$ and $[\mathbf{w}]$ respectively. The distance $\rho_E$ is known as the full *Procrustes distance* (Kendall [22], Kent [23] and Dryden and Mardia [12]).

Let $X_j$, $1 \leq j \leq n$, be i.i.d. k-ads such that their shapes $[X_j]$, $1 \leq j \leq n$, have the common distribution $Q$. Let $\tilde{\mu}$ denote the Euclidean mean of $\tilde{Q} = Q \circ \phi_E^{-1}$ viewed as a probability measure on $S(k, \mathbb{C})$. Let $\tilde{M} = \phi_E(\Sigma_2^k)$, and denote the Euclidean projection of $\tilde{\mu}$ on $\tilde{M}$ by $P\tilde{\mu} \equiv P_{\tilde{M}}\tilde{\mu}$. The *extrinsic mean* of $Q$ is then $\mu_E = \phi^{-1}(P\tilde{\mu})$. It minimizes the *Fréchet function* (2.1) with respect to the distance $\rho_E$. Similarly, for the sample extrinsic mean, calculate $P\overline{\tilde{X}}$ where $\overline{\tilde{X}} = (1/n)\sum_{j=1}^{n} \phi_E([X_j])$ is a coordinate-wise average of the matrix elements $W_j W_j^*$ and $W_j$ is the preshape of $X_j$ ($1 \leq j \leq n$). The asymptotic distribution of $\sqrt{n}(P\overline{\tilde{X}} - P\tilde{\mu})$ is given by that of its projection on the tangent space $T_{P\tilde{\mu}}\tilde{M}$ at $P\tilde{\mu}$, since its projection on the complement of $T_{P\tilde{\mu}}\tilde{M}$ is negligible. For computation of this projection, one chooses a suitable orthonormal basis of $S(k, \mathbb{C})$ (considered as a single orthonormal frame for its constant tangent spaces), and calculates the differential of the projection map $P = P_{\tilde{M}} : S(k, \mathbb{C}) \to \tilde{M}$ in terms of these coordinates. One thus arrives at a nonsingular $(2k - 4)$-dimensional Normal distribution in the limit (see Sections 3.1–3.4 for details).

Turning to the *intrinsic mean* on a Riemannian manifold $M$, with geodesic distance $d_g$, the first problem to resolve is its existence as the unique minimizer of the Fréchet function $\int d_g^2(p, m) Q(dm)$. Here a result of Karchar [19] on the existence of a unique minimizer is greatly improved by a result of Kendall [22], which allows the radius $r$ of a geodesic ball $B(p, r)$ containing the support of $Q$ to be twice as large as required by Karchar [19] (Proposition 4.1). On such a ball, the map $\phi = \exp_p^{-1}$ (the inverse of the *exponential map* at $p$), is a diffeomorphism onto its image in the tangent space $T_p M$ at $p$. Using the coordinates of the vector space $T_p M$, called *normal coordinates*, one arrives at a central limit theorem for the sample intrinsic mean $\mu_{nI}$ (Theorem 4.2), following Bhattacharya and Patrangenaru [8]. Note that, with the (non-Euclidean) distance on $T_p M$ induced by $\phi$ from the geodesic distance $d_g$ on $M$, the image $\mu_n = \phi(\mu_{nI})$ of $\mu_{nI}$ is the minimizer of the Fréchet function

$$F_n(x) \equiv \int d_g^2(\phi^{-1}x, \phi^{-1}y) \tilde{Q}_n(dy)$$

where $\tilde{Q}_n = Q_n \circ \phi^{-1}$, $Q_n = (1/n)\sum_{j=1}^{n} \delta_{[X_j]}$. Thus $\mu_n$ is a *M-estimator* in the Euclidean space $T_q M$. The assumptions in Theorem 4.2 guarantee that this M-estimator is asymptotically Gaussian around $\mu = \phi(\mu_I)$. The asymptotic distribution of the test statistic (4.5) follows from this.

The computation of the test statistic (4.5) is generally more involved than that used for comparing extrinsic means (see, e.g., (3.17) for the case $M = \Sigma_2^k$). This involves, in particular, the metric tensor of $M$ to compute geodesics and normal coordinates. We refer to [3] for the asymptotic theory for intrinsic means, with explicit computations of parameters especially for the planer shape space of k-ads. However in Section 5 of the present article, we display numerical values of the intrinsic two-sample test statistics, along with the corresponding p-values, in two examples. It may be noted that for highly concentrated data in each of these examples, the extrinsic and intrinsic distances are close and hence the extrinsic and intrinsic test statistics have virtually the same values.

The minimum value attained by the Fréchet function is called the *Fréchet variation* of $Q$ and it is a measure of spread of the distribution $Q$. The sample Fréchet variation is a consistent estimator of the Fréchet variation of $Q$ as proved in Proposition 2.4. If the Fréchet mean exists, we derive the asymptotic distribution of the sample Fréchet variation in Theorem 2.5. This can be used to construct a nonparametric test statistic to compare the spread of two populations on $M$. We compute numerical values of the test statistic, along with the p-values for $M = \Sigma_2^k$ in Section 5. For highly concentrated data as in the examples considered in Section 5, the Fréchet variations of the distributions are very small. Then the mean comparison is usually sufficient to discriminate between the populations and the variations show no significant difference.

We conclude this section with two brief remarks. First, the main objective of inference in the two-sample problem on $\Sigma_2^k$ is to discriminate between two different distributions on it. It turns out, in most practical problems that arise, that the means and variations (extrinsic or intrinsic) are generally adequate for this discrimination. More elaborate procedures such as nonparametric density estimation suffer from the "curse of dimensionality" on this commonly high-dimensional space. One can, however, do such density estimation on a tangent space (e.g., on $T_{\mu_I} M$, via the inverse exponential map $\exp_{\mu_I}^{-1}$), as in the Euclidean case. Excepting for the computation in normal coordinates, this presents no novelty. Secondly, in examples with real data sets that we have studied (e.g., those in Section 5), the p-values of the nonparametric two-sample tests for comparing means, developed in this article, are always much smaller (often by an order of magnitude or more) than those based on existing, mostly parametric, tests in the literature (see Dryden and Mardia [12]). This seems to indicate that the tests proposed here may be more powerful than those that have been used in the past, for many data sets that arise in practice. This perhaps also points to the inadequacy of parametric models of shapes popularly used in the literature in capturing certain important shape features.

## 2. Fréchet mean and variation on metric spaces

Let $(M, \rho)$ be a metric space, $\rho$ being the distance on $M$. For a given probability measure $Q$ on (the Borel sigma-field of) $M$, define the *Fréchet function* of $Q$ as

$$(2.1) \qquad F(p) = \int_M \rho^2(p, x) Q(dx), \quad p \in M.$$

### 2.1. Fréchet mean

**Definition 2.1.** Suppose $F(p) < \infty$ for some $p \in M$. Then the set of all $p$ for which $F(p)$ is the minimum value of $F$ on $M$ is called the *Fréchet mean set* of $Q$, denoted by $C_Q$. If this set is a singleton, say $\{\mu_F\}$, then $\mu_F$ is called the *Fréchet mean* of $Q$. If $X_1, X_2, \ldots, X_n$ are independent and identically distributed (i.i.d.) with common distribution $Q$, and $Q_n \doteq (1/n) \sum_{j=1}^n \delta_{X_j}$ is the corresponding empirical distribution, then the Fréchet mean set of $Q_n$ is called the *sample Fréchet mean set*, denoted by $C_{Q_n}$. If this set is a singleton, say $\{\mu_{F_n}\}$, then $\mu_{F_n}$ is called the *sample Fréchet mean*.

The following result has been proved in Theorem 2.1, Bhattacharya and Patrangenaru [7].

**Proposition 2.1.** *Suppose every closed and bounded subset of $M$ is compact. If the Fréchet function $F(p)$ of $Q$ is finite for some $p$, then $C_Q$ is nonempty and compact.*

The next result establishes the strong consistency of the sample Fréchet mean. For a proof, see Theorem 2.3, Bhattacharya and Patrangenaru [7].

**Proposition 2.2.** *Assume* (i) *that every closed bounded subset of $M$ is compact, and* (ii) *$F$ is finite on $M$. Then given any $\epsilon > 0$, there exists an integer valued random variable $N = N(\omega, \epsilon)$ and a $P$-null set $A(\omega, \epsilon)$ such that*

$$(2.2) \qquad C_{Q_n} \subset C_Q^\epsilon \equiv \{p \in M : \rho(p, C_Q) < \epsilon\}, \ \forall n \geq N$$

*outside of $A(\omega, \epsilon)$. In particular, if $C_Q = \{\mu_F\}$, then every measurable selection $\mu_{F_n}$ from $C_{Q_n}$ is a strongly consistent estimator of $\mu_F$.*

**Remark 2.1.** It is known that a connected Riemannian manifold $M$ which is complete (in its geodesic distance) satisfies the topological hypothesis of Propositions 2.1 and 2.2: every closed bounded subset of $M$ is compact (see Theorem 2.8, Do Carmo [10], pp. 146–147). We will investigate conditions for the existence of the Fréchet mean of $Q$ (as a unique minimizer of the Fréchet function $F$ of $Q$) in the subsequent sections.

**Remark 2.2.** One can show that the reverse of (2.2), that is, "$C_Q \subset C_{Q_n}^\epsilon \ \forall \ n \geq N(\omega, \epsilon)$" does not hold in general. See, for example, Bhattacharya and Patrangenaru ([7], Remark 2.6).

Next we consider the asymptotic distribution of $\mu_{F_n}$. For Theorem 2.3, we assume $M$ to be a differentiable manifold of dimension $d$. Let $\rho$ be a distance metrizing the topology of $M$. For a proof of the following result, see Theorem 2.1, Bhattacharya and Patrangenaru [8].

**Theorem 2.3.** *Suppose the following assumptions hold:*

(i) *$Q$ has support in a single coordinate patch, $(U, \phi)$, $\phi : U \longrightarrow \Re^d$ smooth. Let $Y_j = \phi(X_j)$, $j = 1, \ldots, n$.*
(ii) *The Fréchet mean $\mu_F$ of $Q$ is unique.*
(iii) *$\forall x$, $y \mapsto h(x, y) = \rho^2(\phi^{-1}x, \phi^{-1}y)$ is twice continuously differentiable in a neighborhood of $\phi(\mu_F) = \mu$.*
(iv) *$E(\mathrm{D}_r h(Y_1, \mu))^2 < \infty \ \forall r$.*
(v) *$E(\sup_{|u-v|\leq \epsilon} |\mathrm{D}_s \mathrm{D}_r h(Y_1, v) - \mathrm{D}_s \mathrm{D}_r h(Y_1, u)|) \to 0$ as $\epsilon \to 0 \ \forall \ r, s$.*
(vi) *$\Lambda = E(\mathrm{D}_s \mathrm{D}_r h(Y_1, \mu))$ is nonsingular.*
(vii) *$\Sigma = \mathrm{Cov}(\mathrm{D}h(Y_1, \mu))$ is nonsingular.*

*Let $\mu_{F_n}$ be a measurable selection from the Fréchet sample mean set, and write $\mu_n = \phi(\mu_{F_n})$. Then under the assumptions (i)–(vii),*

$$(2.3) \qquad \sqrt{n}(\mu_n - \mu) \xrightarrow{\mathcal{L}} N(0, \Lambda^{-1} \Sigma (\Lambda')^{-1}).$$

### 2.2. Fréchet variation

**Definition 2.2.** The *Fréchet variation $V$* of $Q$ is the minimum value attained by the Fréchet function $F$ defined by (2.1) on $M$. Similarly the minimum value attained by the *sample Fréchet function*,

$$(2.4) \qquad F_n(p) = \frac{1}{n} \sum_{j=1}^n \rho^2(X_j, p)$$

is called the *sample Fréchet variation* and denoted by $V_n$.

From Proposition 2.1 it follows that if the Fréchet function $F(p)$ is finite for some $p$, then $V$ is finite and equals $F(p)$ for all $p$ in the Fréchet mean set $C_Q$. Similarly the sample variation $V_n$ is the value of $F_n$ on the sample Fréchet mean set $C_{Q_n}$. The following result establishes the strong consistency of $V_n$ as an estimator of $V$.

**Proposition 2.4.** *Suppose every closed and bounded subset of $M$ is compact, and $F$ is finite on $M$. Then $V_n$ is a strongly consistent estimator of $V$.*

*Proof.* In view of Proposition 2.2, for any $\epsilon > 0$, there exists $N = N(\omega, \epsilon)$ such that

$$(2.5) \qquad |V_n - V| = |\inf_{p \in M} F_n(p) - \inf_{p \in M} F(p)| \leq \sup_{p \in \overline{C_Q^\epsilon}} |F_n(p) - F(p)|$$

for all $n \geq N$ almost surely. From the proof of Theorem 2.3 in Bhattacharya and Patrangenaru [7], it follows that for any compact set $K \subset M$,

$$\sup_{p \in K} |F_n(p) - F(p)| \longrightarrow 0 \text{ a.s. as } n \to \infty.$$

Since $\overline{C_Q^\epsilon}$ is compact, it follows from (2.5) that

$$|V_n - V| \longrightarrow 0 \text{ a.s. as } n \to \infty. \qquad \square$$

**Remark 2.3.** The sample variation is a consistent estimator of the population variation even when the Fréchet function $F$ of $Q$ does not have a unique minimizer.

Next we derive the asymptotic distribution of $V_n$ when there is a unique population Fréchet mean.

**Theorem 2.5.** *Let $M$ be a differentiable manifold. Using the notation of Theorem 2.3, under assumptions (i)–(vii) and assuming $E(\rho^4(X_1, \mu_F)) < \infty$, one has*

$$(2.6) \qquad \sqrt{n}(V_n - V) \xrightarrow{\mathcal{L}} N\left(0, \mathrm{var}(\rho^2(X_1, \mu_F))\right).$$

*Proof.* Let

$$F(x) = \int \rho^2(\phi^{-1}(x), m) Q(dm), \quad F_n(x) = \frac{1}{n} \sum_{j=1}^n \rho^2(\phi^{-1}(x), X_j).$$

Let $\mu_{F_n}$ be a measurable selection from the sample mean set and $\mu_n = \phi(\mu_{Fn})$. Then

$$\sqrt{n}(V_n - V) = \sqrt{n}(F_n(\mu_n) - F(\mu))$$
$$(2.7) \qquad = \sqrt{n}(F_n(\mu_n) - F_n(\mu)) + \sqrt{n}(F_n(\mu) - F(\mu)),$$

$$\sqrt{n}(F_n(\mu_n) - F_n(\mu)) = \frac{1}{\sqrt{n}} \sum_{j=1}^n \sum_{r=1}^d (\mu_n - \mu)_r D_r h(Y_j, \mu)$$
$$(2.8) \qquad + \frac{1}{2\sqrt{n}} \sum_{j=1}^n \sum_{r=1}^d \sum_{s=1}^d (\mu_n - \mu)_r (\mu_n - \mu)_s D_s D_r h(Y_j, \mu_n^*)$$

for some $\mu_n^*$ in the line segment joining $\mu$ and $\mu_n$. By assumption (v) of Theorem 2.3 and because $\sqrt{n}(\mu_n - \mu)$ is asymptotically normal, the second term on the right of (2.8) converges to 0 in probability. Also $(1/n)\sum_{j=1}^{n} Dh(Y_j, \mu) \to E(Dh(Y_1, \mu)) = 0$, so that the first term on the right of (2.8) converges to 0 in probability. Hence (2.7) becomes

$$\sqrt{n}(V_n - V) = \sqrt{n}(F_n(\mu) - F(\mu)) + o_P(1)$$
$$(2.9) \qquad = \frac{1}{\sqrt{n}} \sum_{j=1}^{n} \left(\rho^2(X_j, \mu_F) - E\rho^2(X_1, \mu_F)\right) + o_P(1).$$

By the CLT for the i.i.d. sequence $\{\rho^2(X_j, \mu_F)\}$, (2.9) converges in law to $N(0, \text{var}(\rho^2(X_1, \mu_F)))$. □

**Remark 2.4.** Theorem 2.5 requires the population mean to exist for the sample variation to be asymptotically Normal. It may be shown by examples that it fails to give the correct distribution if there is not a unique mean.

Theorem 2.5 can be used to construct a nonparametric test for testing whether two populations have the same spread. Suppose $Q_1$ and $Q_2$ are two probability distributions with unique Fréchet means $\mu_{1F}$ and $\mu_{2F}$ and Fréchet variations $V_1$ and $V_2$, respectively. We have i.i.d. samples $X_1, X_2, \ldots, X_n$ and $Y_1, Y_2, \ldots, Y_m$ from $Q_1$ and $Q_2$, respectively. Let $\mu_{F_n}$ and $\mu_{F_m}$ denote the sample means, $V_n$ and $V_m$ denote the sample variations. Then the null hypothesis is

$$H_0 : V_1 = V_2 = V.$$

Under $H_0$, from (2.6),

$$(2.10) \qquad \sqrt{n}(V_n - V) \xrightarrow{\mathcal{L}} N(0, \sigma_1^2),$$

$$(2.11) \qquad \sqrt{m}(V_m - V) \xrightarrow{\mathcal{L}} N(0, \sigma_2^2),$$

where $\sigma_1^2 = \text{var}(\rho^2(X_1, \mu_{1F}))$, $\sigma_2^2 = \text{var}(\rho^2(Y_1, \mu_{2F}))$.

Suppose $n/(m+n) \to p$, $m/(m+n) \to q$, for some $p, q > 0$; $p + q = 1$. Then from (2.10) and (2.11),

$$(2.12) \qquad \sqrt{n+m}(V_n - V_m) \xrightarrow{\mathcal{L}} N\left(0, \left(\frac{\sigma_1^2}{p} + \frac{\sigma_2^2}{q}\right)\right),$$

$$(2.13) \qquad \frac{V_n - V_m}{\sqrt{\frac{s_1^2}{n} + \frac{s_2^2}{m}}} \xrightarrow{\mathcal{L}} N(0, 1),$$

where $s_1^2 = (1/n)\sum_{j=1}^{n}(\rho^2(X_j, \mu_{F_n}) - V_n)^2$ and $s_2^2 = (1/m)\sum_{j=1}^{m}(\rho^2(Y_j, \mu_{F_m}) - V_m)^2$ are the sample estimates of $\sigma_1^2$ and $\sigma_2^2$, respectively. Hence the test statistic used is

$$(2.14) \qquad T_{nm} = \frac{V_n - V_m}{\sqrt{\frac{s_1^2}{n} + \frac{s_2^2}{m}}}.$$

For a test of size $\alpha$, we reject $H_0$ if $|T_{nm}| > Z_{1-(\alpha/2)}$ where $Z_{1-(\alpha/2)}$ is the $(1-(\alpha/2))^{\text{th}}$ quantile of $N(0,1)$.

From now on, unless otherwise stated, we assume that $(\mathbf{M}, \mathbf{g})$ is a $d$-dimensional connected complete Riemannian manifold, $g$ being the Riemannian metric tensor on $M$. We shall come across different notions of means and variations depending on the distance chosen on $M$. We begin with the *extrinsic distance* in the next section.

## 3. Extrinsic mean and variation

Let $\phi : M \to \Re^k$ be an embedding of $M$ into $\Re^k$, and let $\tilde{M} = \phi(M) \subset \Re^k$. Define the distance on $M$ as: $\rho(x,y) = \|\phi(x) - \phi(y)\|$, where $\|\cdot\|$ denotes Euclidean norm ($\|u\|^2 = \sum_{i=1}^k u_i^2$, $u = (u_1, u_2, \ldots, u_k)'$). This is called the *extrinsic distance* on $M$.

Assume that $\tilde{M}$ is a closed subset of $\Re^k$. Then for every $u \in \Re^k$ there exists a compact set of points in $\tilde{M}$ whose distance from $u$ is the smallest among all points in $\tilde{M}$. We will denote this set by

$$Pu \equiv P_{\tilde{M}} u = \{x \in \tilde{M} : \|x - u\| \leq \|y - u\| \; \forall y \in \tilde{M}\}.$$

If this set is a singleton, $u$ is said to be a *nonfocal point* of $\Re^k$ (with respect to $\tilde{M}$); otherwise it is said to be a *focal point* of $\Re^k$.

**Definition 3.1.** Let $(M, \rho)$, $\phi$ be as above. Let $Q$ be a probability measure on $M$ with finite Fréchet function. The Fréchet mean (set) of $Q$ is called the *extrinsic mean* (set) of $Q$, and the Fréchet variation of $Q$ is called its *extrinsic variation*. If $X_j$ ($j = 1, \ldots, n$) are iid observations from $Q$, and $Q_n = \frac{1}{n} \sum_{j=1}^n \delta_{X_j}$ is the empirical distribution, then the Fréchet mean(set) of $Q_n$ is called the *extrinsic sample mean*(set) and the Fréchet variation of $Q_n$ is called the *extrinsic sample variation*.

Let $\tilde{Q}$ and $\tilde{Q}_n$ be the images of $Q$ and $Q_n$, respectively, on $\Re^k$ under $\phi$: $\tilde{Q} = Q \circ \phi^{-1}$, $\tilde{Q}_n = Q_n \circ \phi^{-1}$. The next result gives us a way to calculate the extrinsic mean and establishes the consistency of the sample mean as an estimator of the population mean if that exists. For a proof see Proposition 3.1 in Bhattacharya and Patrangenaru [7].

**Proposition 3.1.** (a) If $\tilde{\mu} = \int_{R^k} u \tilde{Q}(du)$ is the mean of $\tilde{Q}$, then the extrinsic mean set of $Q$ is given by $\phi^{-1}(P\tilde{\mu})$. (b) If $\tilde{\mu}$ is a nonfocal point of $\Re^k$ (relative to $\tilde{M}$), then the extrinsic sample mean $\mu_{nE}$ (any measurable selection from the extrinsic mean set of $Q_n$) is a strongly consistent estimator of the extrinsic mean $\mu_E = \phi^{-1}(P\tilde{\mu})$.

### 3.1. Asymptotic distribution of the sample extrinsic mean

We can use Theorem 2.3 to get the asymptotic distribution of the sample extrinsic mean. However, expressions for the parameters $\Lambda$ and $\Sigma$ are not easy to get. Here we devise another way to derive the asymptotic distribution. We assume that the mean $\tilde{\mu}$ of $\tilde{Q}$ is a nonfocal point, so that the projection $P\tilde{\mu}$ of $\tilde{\mu}$ on $\phi(M)$ is unique, and the extrinsic mean of $Q$ is $\mu_E = \phi^{-1}(P\tilde{\mu})$. Let $\overline{\tilde{X}} = (1/n) \sum_{j=1}^n \tilde{X}_j$ denote the sample mean of $\tilde{X}_j = \phi(X_j)$. The extrinsic sample mean set is $C_{Q_n} = \phi^{-1}(P\overline{\tilde{X}})$, where $P\overline{\tilde{X}}$ is the set of projection of $\overline{\tilde{X}}$ on $\phi(M)$. In a neighborhood of a nonfocal point such as $\tilde{\mu}$, $P(.)$ is smooth. So we can write

$$(3.1) \quad \sqrt{n}[P(\overline{\tilde{X}}) - P(\tilde{\mu})] = \sqrt{n}(d_{\tilde{\mu}} P)(\overline{\tilde{X}} - \tilde{\mu}) + o_P(1) = (d_{\tilde{\mu}} P)(\sqrt{n}(\overline{\tilde{X}} - \tilde{\mu})) + o_P(1)$$

where $d_{\tilde{\mu}} P$ is the differential (map) of $P(\cdot)$, which takes vectors in the tangent space of $\Re^k$ at $\tilde{\mu}$ to tangent vectors of $\phi(M)$ at $P(\tilde{\mu})$. Hence the left side is asymptotically normal.

For the case of regular submanifolds embedded in an Euclidean space by the inclusion map, a similar asymptotic distribution and a two-sample test were constructed independently by Hendricks and Landsman [18] and, for more general manifolds, by Patrangenaru [26] and Bhattacharya and Patrangenaru [8].

### 3.2. Application to the planar shape space of k-ads

Consider a set of $k$ points on the plane, e.g., $k$ locations on a skull projected on a plane, not all points being the same. We will assume $k > 2$ and refer to such a set as a *k-ad* (or a set of *k landmarks*). For convenience we will denote a k-ad by $k$ complex numbers $(z_j = x_j + i y_j, 1 \leq j \leq k)$, i.e., we will represent k-ads on a complex plane. By the *shape* of a k-ad $\mathbf{z} = (z_1, z_2, \ldots, z_k)$, we mean the equivalence class, or orbit of $\mathbf{z}$ under translation, rotation and scaling. To remove translation, one may substract $\langle \mathbf{z} \rangle \equiv (\langle z \rangle, \langle z \rangle, \ldots, \langle z \rangle)$ $(\langle z \rangle = (1/k) \sum_{j=1}^{k} z_j)$ from $\mathbf{z}$ to get $\mathbf{z} - \langle \mathbf{z} \rangle$. Rotation of the k-ad by an angle $\theta$ and scaling (by a factor $r > 0$) are achieved by multiplying $\mathbf{z} - \langle \mathbf{z} \rangle$ by the complex number $\lambda = r \exp i\theta$. Hence one may represent the shape of the k-ad as the complex line passing through $\mathbf{z} - \langle \mathbf{z} \rangle$, namely, $\{\lambda(\mathbf{z} - \langle \mathbf{z} \rangle) \colon \lambda \in \mathbb{C} \setminus \{0\}\}$. Thus the space of k-ads is the set of all complex lines on the (complex $(k-1)$-dimensional) hyperplane, $H^{k-1} = \{w \in C^k \setminus \{0\} \colon \sum_{1}^{k} w_j = 0\}$. Therefore the shape space $\Sigma_2^k$ of planer k-ads has the structure of the *complex projective space* $\mathbb{C}P^{k-2}$: the space of all complex lines through the origin in $\mathbb{C}^{k-1}$. As in the case of $\mathbb{C}P^{k-2}$, it is convenient to represent the element of $\Sigma_2^k$ corresponding to a k-ad $\mathbf{z}$ by the curve $\gamma(z) = [z] = \{e^{i\theta}((z - \langle \mathbf{z} \rangle)/\|z - \langle \mathbf{z} \rangle\|) \colon 0 \leq \theta < 2\pi\}$ on the unit sphere in $H^{k-1} \approx \mathbb{C}^{k-1}$.

If we denote by $u$ the quantity $(\mathbf{z} - \langle \mathbf{z} \rangle)/\|\mathbf{z} - \langle \mathbf{z} \rangle\|$, called the *preshape* of the shape of $\mathbf{z}$, then another representation of $\Sigma_2^k$ is via the *Veronese–Whitney embedding* $\phi$ into the space $S(k, \mathbb{C})$ of all $k \times k$ complex Hermitian matrices. $S(k, \mathbb{C})$ is viewed as a (real) vector space with respect to the scaler field $\Re$. The embedding $\phi$ is given by

$$\phi \colon \Sigma_2^k \to S(k, \mathbb{C}),$$
$$\phi([z]) = uu^* \ (u = (u_1, \ldots, u_k)' \in H^{k-1}, \|u\| = 1)$$
(3.2)
$$= ((u_i \bar{u}_j))_{1 \leq i, j \leq k}.$$

The shape of $\mathbf{z}$, $[z] = \{e^{i\theta} u \colon 0 \leq \theta < 2\pi\}$ is the orbit of the vector $u$ under rotation. Note that if $v_1, v_2 \in [z]$, then $\phi([v_1]) = \phi([v_2]) = \phi((z - \langle \mathbf{z} \rangle)/\|z - \langle \mathbf{z} \rangle\|)$. Define the *extrinsic distance* $\rho$ on $\Sigma_2^k$ by that induced from this embedding, namely,

$$(3.3) \qquad \rho^2([z], [w]) = \|uu^* - vv^*\|^2 \ , u \doteq \frac{z - \langle \mathbf{z} \rangle}{\|z - \langle \mathbf{z} \rangle\|} \ , v \doteq \frac{w - \langle \mathbf{w} \rangle}{\|w - \langle \mathbf{w} \rangle\|}$$

where for arbitrary $k \times k$ complex matrices A, B,

$$(3.4) \qquad \|A - B\|^2 = \sum_{j,j'} |a_{jj'} - b_{jj'}|^2 = \text{Trace}(A - B)(A - B)^*$$

is just the squared euclidean distance between A and B regarded as elements of $\mathbb{C}^{k^2}$ (or, $\Re^{2k^2}$). Since the matrices $uu^*$, $vv^*$ in (3.2) are Hermitian, one notes that the image $\phi(\Sigma_2^k)$ of $\Sigma_2^k$ is a closed subset of $\mathbb{C}^{k^2}$ and the "conjugate-transpose" symbol $*$ may be dropped from (3.4) in computing distances in $\phi(\Sigma_2^k)$.

Let $Q$ be a probability measure on the shape space $\Sigma_2^k$, let $[X_1], [X_2], \ldots, [X_n]$ be an i.i.d. sample from $Q$ and let $\tilde{\mu}$ denote the mean vector of $\tilde{Q} \doteq Q \circ \phi^{-1}$, regarded as a probability measure on $\mathbb{C}^{k^2}$ (or, $\mathfrak{R}^{2k^2}$). Note that $\tilde{\mu}$ belongs to the convex hull of $\tilde{M} = \phi(\Sigma_2^k)$ and in particular, is an element of $H^{k-1}$. Let $T$ be a (complex) orthogonal $k \times k$ matrix such that $T\tilde{\mu}T^* = D = \text{Diag}(\lambda_1, \lambda_2, \ldots, \lambda_k)$, where $\lambda_1 \leq \lambda_2 \leq \cdots \leq \lambda_k$ are the eigenvalues of $\tilde{\mu}$. Then, writing $v = Tu$ with $u$ as in (3.3),

$$\|uu^* - \tilde{\mu}\|^2 = \|vv^* - D\|^2 = \sum_{j=1}^k (|v_j|^2 - \lambda_j)^2 + \sum_{j \neq j'} |v_j \overline{v}_{j'}|^2$$

$$= \sum \lambda_j^2 + \sum_{j=1}^k |v_j|^4 - 2\sum_{j=1}^k \lambda_j |v_j|^2 + \sum_{j=1}^k |v_j|^2 \cdot \sum_{j'=1}^k |v_{j'}|^2 - \sum_{j=1}^k |v_j|^4$$

(3.5)
$$= \sum \lambda_j^2 + 1 - 2\sum_{j=1}^k \lambda_j |v_j|^2$$

which is minimized (on $\phi(\Sigma_2^k)$) by taking $v = e_k = (0, \ldots, 0, 1)'$, i.e., $u = T^* e_k$, a unit eigenvector having the largest eigenvalue $\lambda_k$ of $\tilde{\mu}$. It follows that the extrinsic mean $\mu_E$, say, of $Q$ is unique if and only if the eigenspace for the largest eigenvalue of $\tilde{\mu}$ is (complex) one-dimensional, and then $\mu_E = [\mu]$, $\mu (\neq 0) \in$ the eigenspace of the largest eigenvalue of $\tilde{\mu}$.

From (3.5), the extrinsic variation of $Q$ has the expression

$$V = \mathrm{E}\|X_1 X_1^* - \mu\mu^*\|^2$$
$$= \mathrm{E}\|X_1 X_1^* - \tilde{\mu}\|^2 + \|\tilde{\mu} - \mu\mu^*\|^2$$
(3.6)
$$= 2(1 - \lambda_k)$$

Therefore, we have the following consequence of Proposition 2.4 and Proposition 3.1.

**Corollary 3.2.** *Let $\mu_n$ denote an eigenvector of $(1/n)\sum_{j=1}^n X_j X_j^*$ having the largest eigenvalue $\lambda_{kn}$. (a) If the largest eigenvalue $\lambda_k$ of $\tilde{\mu}$ is simple, then the extrinsic sample mean $[\mu_n]$ is a strongly consistent estimator of the extrinsic mean $[\mu]$. (b) The sample extrinsic variation, $V_n = 2(1 - \lambda_{kn})$ is a strongly consistent estimator of the extrinsic variation, $V = 2(1 - \lambda_k)$.*

The distance $\rho$ on $\Sigma_2^k$ in (3.3) can be expressed as

(3.7) $$\rho^2([z], [w]) \equiv \|uu^* - vv^*\|^2 = 2(1 - |u^* v|^2).$$

This is the so-called *full Procrustes distance* for $\Sigma_2^k$. See Kent [23], Dryden and Mardia [12] and Kendall et al. [21].

### 3.3. Asymptotic distribution of mean shape

To get the asymptotic distribution of the sample extrinsic mean shape using (3.1), we embed $M = \Sigma_2^k$ into $S(k, \mathbb{C})$, the space of all $k \times k$ complex self-adjoint matrices, via the map $\phi$ in (3.2). We consider $S(k, \mathbb{C})$ as a linear subspace of $\mathbb{C}^{k^2}$ (over $\mathfrak{R}$) and as such a regular submanifold of $\mathbb{C}^{k^2}$ embedded by the inclusion map, and inheriting the metric tensor:

$$\langle A, B \rangle = \text{Re}\left(\text{Trace}(A\bar{B}')\right).$$

The (real) dimension of $S(k, \mathbb{C})$ is $k^2$. An orthonormal basis for $S(k, \mathbb{C})$ is given by $\{v_b^a : 1 \leq a \leq b \leq k\}$ and $\{w_b^a : 1 \leq a < b \leq k\}$, defined as

$$v_b^a = \begin{cases} \frac{1}{\sqrt{2}}(e_a e_b^t + e_b e_a^t), & a < b \\ e_a e_a^t, & a = b \end{cases}$$

$$w_b^a = +\frac{i}{\sqrt{2}}(e_a e_b^t - e_b e_a^t),\ a < b.$$

where $\{e_a : 1 \leq a \leq k\}$ is the standard canonical basis for $\Re^k$.

We also take $\{v_b^a : 1 \leq a \leq b \leq k\}$ and $\{w_b^a : 1 \leq a < b \leq k\}$ as the (constant) orthogonal frame for $S(k, \mathbb{C})$. For any $U \in O(k)$ ($UU^* = U^*U = I$), $\{Uv_b^a U^* : 1 \leq a \leq b \leq k\}$, $\{Uw_b^a U^* : 1 \leq a < b \leq k\}$ is also an orthogonal frame for $S(k, \mathbb{C})$. Assume that the mean $\tilde{\mu}$ of $\tilde{Q}$ has its largest eigenvalue simple. To apply (3.1), we view $d_{\tilde{\mu}} P : S(k, \mathbb{C}) \to T_{P(\tilde{\mu})} \phi(\Sigma_2^k)$. Choose $U \in O(k)$ such that $U^* \tilde{\mu} U = D \equiv \text{Diag}(\lambda_1, \ldots, \lambda_k)$, $\lambda_1 \leq \cdots \leq \lambda_{k-1} < \lambda_k$ being the eigenvalues of $\tilde{\mu}$.

Choose the basis frame $\{Uv_b^a U^*, Uw_b^a U^*\}$ for $S(k, \mathbb{C})$. Then one can show that

$$d_{\tilde{\mu}} P(Uv_b^a U^*) = \begin{cases} 0, & \text{if } 1 \leq a \leq b < k,\ a = b = k, \\ (\lambda_k - \lambda_a)^{-1} Uv_k^a U^*, & \text{if } 1 \leq a < k, b = k. \end{cases}$$

(3.8) $$d_{\tilde{\mu}} P(Uw_b^a U^*) = \begin{cases} 0, & \text{if } 1 \leq a < b < k, \\ (\lambda_k - \lambda_a)^{-1} Uw_k^a U^*, & \text{if } 1 \leq a < k, b = k. \end{cases}$$

Write

$$\sqrt{n}(\bar{\tilde{X}} - \tilde{\mu}) = \sum_{1 \leq a \leq b \leq k} \langle \sqrt{n}(\bar{\tilde{X}} - \tilde{\mu}), Uv_b^a U^* \rangle Uv_b^a U^*$$

(3.9) $$+ \sum_{1 \leq a < b \leq k} \langle \sqrt{n}(\bar{\tilde{X}} - \tilde{\mu}), Uw_b^a U^* \rangle Uw_b^a U^*.$$

Since $\bar{\tilde{X}} \mathbf{1}_k = \tilde{\mu} \mathbf{1}_k = 0$, $\lambda_1 = 0$ and $U_{\cdot 1} = \alpha \mathbf{1}_k$, $|\alpha| = 1/\sqrt{k}$. Thus

$$\langle \sqrt{n}(\bar{\tilde{X}} - \tilde{\mu}), Uv_b^1 U^* \rangle = \langle \sqrt{n}(\bar{\tilde{X}} - \tilde{\mu}), Uw_b^1 U^* \rangle = 0.$$

Therefore,

$$d_{\tilde{\mu}} P(\sqrt{n}(\bar{\tilde{X}} - \tilde{\mu}))$$
$$= \sum_{a=2}^{k-1} \langle \sqrt{n}(\bar{\tilde{X}} - \tilde{\mu}), Uv_k^a U^* \rangle (\lambda_k - \lambda_a)^{-1} Uv_k^a U^*$$

(3.10) $$+ \sum_{a=2}^{k-1} \langle \sqrt{n}(\bar{\tilde{X}} - \tilde{\mu}), Uw_k^a U^* \rangle (\lambda_k - \lambda_a)^{-1} Uw_k^a U^*.$$

From (3.10), we see that $\sqrt{n}(P(\bar{\tilde{X}}) - P(\tilde{\mu}))$ has an asymptotic Gaussian distribution on a subspace of $S(k, \mathbb{C})$ with asymptotic coordinates

$$T_n(\tilde{\mu}) = \left( \langle \sqrt{n}(\bar{\tilde{X}} - \tilde{\mu}), Uv_k^a U^* \rangle_{a=2}^{k-1},\ \langle \sqrt{n}(\bar{\tilde{X}} - \tilde{\mu}), Uw_k^a U^* \rangle_{a=2}^{k-1} \right)$$

with respect to the basis vector $\{(\lambda_k - \lambda_a)^{-1} Uv_k^a U^*, (\lambda_k - \lambda_a)^{-1} Uw_k^a U^*\}_{a=2}^{k-1}$.

Writing $\Sigma(\tilde{\mu})$ for the covariance matrix of $T_n(\tilde{\mu})$, and assuming that it is non-singular,

(3.11) $$T_n(\tilde{\mu})' \Sigma(\tilde{\mu})^{-1} T_n(\tilde{\mu}) \longrightarrow \mathcal{X}_{2k-4}^2.$$

### 3.4. Two sample testing problems on $\Sigma_2^k$

Let $Q_1$ and $Q_2$ be two probability measures on the shape space $\Sigma_2^k$, and let $\mu_1$ and $\mu_2$ denote the means of $Q_1 \circ \phi^{-1}$ and $Q_2 \circ \phi^{-1}$, respectively. Suppose $[x_1], \ldots, [x_n]$ and $[y_1], \ldots, [y_m]$ are i.i.d. random samples from $Q_1$ and $Q_2$ respectively. Let $X_i = \phi([x_i])$, $Y_i = \phi([y_i])$ be their images onto $\phi(\Sigma_2^k)$ which are random samples from $Q_1 \circ \phi^{-1}$ and $Q_2 \circ \phi^{-1}$, respectively. Suppose we are to test if the extrinsic means of $Q_1$ and $Q_2$ are equal, i.e.

$$H_0 : P\mu_1 = P\mu_2$$

We assume that both $\mu_1$ and $\mu_2$ have simple largest eigenvalues. Then under $H_0$, the corresponding eigenvectors differ by a rotation.

Choose $\mu \in S(k, \mathbb{C})$ with the same projection as $\mu_1$ and $\mu_2$. Suppose $\mu = U\Lambda U^*$, where $\Lambda = \mathrm{Diag}(\lambda_1 \leq \lambda_2 \leq \cdots < \lambda_k)$ are its eigenvalues and $U = [U_1, U_2, \ldots, U_k]$ are the corresponding eigenvectors. Under $H_0$, $P\mu_1 = P\mu_2 = U_k U_k^*$. From (3.10),

$$d_\mu P(\bar{X} - \mu)$$

$$= \sum_{a=2}^{k-1} \sqrt{2}\mathrm{Re}(U_a^* \bar{X} U_k)(\lambda_k - \lambda_a)^{-1} U v_k^a U^*$$

$$+ \sum_{a=2}^{k-1} \sqrt{2}\mathrm{Im}(U_a^* \bar{X} U_k)(\lambda_k - \lambda_a)^{-1} U w_k^a U^*$$

$$= \sum_{a=2}^{k-1} (\lambda_k - \lambda_a)^{-1}(U_a^* \bar{X} U_k) U_a U_k^*$$

(3.12)
$$+ \sum_{a=2}^{k-1} (\lambda_k - \lambda_a)^{-1}(U_k^* \bar{X} U_a) U_k U_a^*,$$

$$d_\mu P(\bar{Y} - \mu)$$

$$= \sum_{a=2}^{k-1} \sqrt{2}\mathrm{Re}(U_a^* \bar{Y} U_k)(\lambda_k - \lambda_a)^{-1} U v_k^a U^*$$

$$+ \sum_{a=2}^{k-1} \sqrt{2}\mathrm{Im}(U_a^* \bar{Y} U_k)(\lambda_k - \lambda_a)^{-1} U w_k^a U^*$$

$$= \sum_{a=2}^{k-1} (\lambda_k - \lambda_a)^{-1}(U_a^* \bar{Y} U_k) U_a U_k^*$$

(3.13)
$$+ \sum_{a=2}^{k-1} (\lambda_k - \lambda_a)^{-1}(U_k^* \bar{Y} U_a) U_k U_a^*.$$

Define

$$T(\mu)_{ij} = \begin{cases} \mathrm{Re}(U_{i+1}^* X_j U_k), & \text{if } 1 \leq i \leq k-2,\ 1 \leq j \leq n, \\ \mathrm{Im}(U_{i-k+3}^* X_j U_k), & \text{if } k-1 \leq i \leq 2k-4,\ 1 \leq j \leq n, \end{cases}$$

$$S(\mu)_{ij} = \begin{cases} \mathrm{Re}(U_{i+1}^* Y_j U_k), & \text{if } 1 \leq i \leq k-2,\ 1 \leq j \leq m, \\ \mathrm{Im}(U_{i-k+3}^* Y_j U_k), & \text{if } k-1 \leq i \leq 2k-4,\ 1 \leq j \leq m, \end{cases}$$

(3.14) $\quad \bar{T}(\mu) = \dfrac{1}{n} \sum_{j=1}^{n} T(\mu)_{.j},\ \bar{S}(\mu) = \dfrac{1}{m} \sum_{j=1}^{m} S(\mu)_{.j}.$

Under $H_0$, $\bar{T}(\mu)$ and $\bar{S}(\mu)$ have mean zero, and as $n, m \to \infty$,

$$\text{(3.15)} \quad \sqrt{n}\bar{T}(\mu) \xrightarrow{\mathcal{L}} N(0, \Sigma_1(\mu)), \; \sqrt{m}\bar{S}(\mu) \xrightarrow{\mathcal{L}} N(0, \Sigma_2(\mu))$$

where $\Sigma_1(\mu)$ and $\Sigma_2(\mu)$ are the covariances of $T(\mu)_{.1}$ and $S(\mu)_{.1}$, respectively. Suppose $(n/(m+n)) \to p$, $(m/(m+n)) \to q$, for some $p, q > 0$; $p + q = 1$. Then

$$\sqrt{n+m}(\bar{T}(\mu) - \bar{S}(\mu)) \xrightarrow{\mathcal{L}} N_{2k-4}(0, \frac{1}{p}\Sigma_1(\mu) + \frac{1}{q}\Sigma_2(\mu)).$$

Thus assuming $\Sigma_1(\mu), \Sigma_2(\mu)$ and hence $\frac{1}{p}\Sigma_1(\mu) + \frac{1}{q}\Sigma_2(\mu)$ to be nonsingular,

$$\text{(3.16)} \quad (n+m)(\bar{T}(\mu) - \bar{S}(\mu))'(\frac{1}{p}\Sigma_1(\mu) + \frac{1}{q}\Sigma_2(\mu))^{-1}(\bar{T}(\mu) - \bar{S}(\mu)) \xrightarrow{\mathcal{L}} \mathcal{X}^2_{2k-4}.$$

Note that the nonsingularity assumption for $\Sigma_1(\mu)$ and $\Sigma_2(\mu)$ are satisfied if, for example, $Q_1$ and $Q_2$ have nonzero absolutely continuous components with respect to the volume measure on $\Sigma_2^k$ (identified with the Riemannian manifold $\mathbb{C}P^{k-2}$). We can choose $\mu$ to be any positive linear combination of $\mu_1$ and $\mu_2$. Then under $H_0$, $\mu$ will have the same projection on $\phi(\Sigma_2^k)$ as $\mu_1$ and $\mu_2$. We may take $\mu = p\mu_1 + q\mu_2$. In practice, since $\mu_1$ and $\mu_2$ are unknown, so is $\mu$. Then we may estimate $\mu$ by the pooled sample mean $\hat{\mu} = (n\bar{X} + m\bar{Y})/(m+n)$, $\Sigma_1(\mu)$ and $\Sigma_2(\mu)$ by their sample estimates $\hat{\Sigma}_1(\hat{\mu})$ and $\hat{\Sigma}_2(\hat{\mu})$, where

$$\hat{\Sigma}_1(\mu) = \frac{1}{n}T(\mu)T(\mu)' - \bar{T}(\mu)\bar{T}(\mu)', \; \hat{\Sigma}_2(\mu) = \frac{1}{m}S(\mu)S(\mu)' - \bar{S}(\mu)\bar{S}(\mu)'.$$

Then the two-sample test statistic in (3.16) can be estimated by

$$\text{(3.17)} \quad T_{nm} = (\bar{T}(\hat{\mu}) - \bar{S}(\hat{\mu}))'(\frac{1}{n}\hat{\Sigma}_1(\hat{\mu}) + \frac{1}{m}\hat{\Sigma}_2(\hat{\mu}))^{-1}(\bar{T}(\hat{\mu}) - \bar{S}(\hat{\mu})).$$

Given level $\alpha$, we reject $H_0$ if

$$\text{(3.18)} \quad T_{nm} > \mathcal{X}^2_{2k-4}(1-\alpha).$$

The expression for $T_{nm}$ depends on the spectrum of $\hat{\mu}$ through the orbit $[U_k(\hat{\mu})]$ and the subspace spanned by $\{U_2(\hat{\mu}), \ldots, U_{k-1}(\hat{\mu})\}$. If the population mean exists, $[U_k(\hat{\mu})]$ is a consistent estimator of $[U_k(\mu)]$ and by perturbation theory (see Dunford and Schwartz [13], p. 598), the projection on $\text{Span}\{U_2(\hat{\mu}), \ldots, U_{k-1}(\hat{\mu})\}$ converges to that on $\text{Span}\{U_2(\mu), \ldots, U_{k-1}(\mu)\}$. Thus from (3.16) and (3.17), $T_{nm}$ has an asymptotic $\mathcal{X}^2_{2k-4}$ distribution. Hence the test in (3.18) has asymptotic level $\alpha$.

To test if the populations have the same spread around their respective means, we use the test statistic in (2.14), which is

$$\text{(3.19)} \quad T_{nm} = 2\frac{\lambda_{km} - \lambda_{kn}}{\sqrt{\frac{s_1^2}{n} + \frac{s_2^2}{m}}},$$

where $\lambda_{kn}$ and $\lambda_{km}$ are the largest eigenvalues of $\bar{X}$ and $\bar{Y}$, respectively. Under $H_0$, $T_{nm}$ has asymptotic Normal distribution.

## 4. Intrinsic mean and variation

Let $(M, g)$ be a d-dimensional connected complete Riemannian manifold, $g$ being the Riemannian metric on $M$. Let the distance $\rho = d_g$ be the geodesic distance under $g$. Let $Q$ be a probability distribution on $M$ with finite Fréchet function,

$$F(p) = \int_M d_g^2(p, m) Q(dm), \; p \in M. \tag{4.1}$$

**Definition 4.1.** The Fréchet mean (set) of $Q$ under the distance $d_g$ is called its *intrinsic mean* (set). The Fréchet variation of $Q$ under $d_g$ is called its *intrinsic variation*. Let $X_1, X_2, \ldots, X_n$ be i.i.d. observations on $M$ with common distribution $Q$. The sample Fréchet mean (set) is called the *sample intrinsic mean* (set) and the sample Fréchet variation is called the *sample intrinsic variation*.

Let us define a few technical terms related to Riemannian manifolds which we will use extensively in the subsequent sections. For details on Riemannian Manifolds, see DoCarmo [10], Gallot et al. [15] or Lee [25].

1. *Geodesic*: These are curves $\gamma$ on the manifold with zero acceleration. They are locally length minimizing curves. For example, consider great circles on the sphere or straight lines in $\Re^d$.
2. *Exponential map*: For $p \in M$, $v \in T_p M$, we define $\exp_p v = \gamma(1)$, where $\gamma$ is a geodesic with $\gamma(0) = p$ and $\dot\gamma(0) = v$.
3. *Cut locus*: Let $\gamma$ be a unit speed geodesic starting at p, $\gamma(0) = p$. Let $t_0$ be the supremum of all $t$ for which $\gamma$ is length minimizing on $[0, t]$. Then $\gamma(t_0)$ is called the cut point of $p$ along $\gamma$. The *cut locus* of $p$, $C(p)$, is the set of all cut points of $p$ along all geodesics. For example, $C(p) = \{-p\}$ on $S^d$.
4. *Convex ball*: A ball $B$ is called convex if, for any $p, q \in B$, a unique geodesic from $p$ to $q$ lies in $B$, which is also the shortest geodesic from $p$ to $q$. For example, any ball of radius $\pi/2$ or less in $S^d$ is convex.
5. *Sectional Curvature*: Recall the notion of Gaussian curvature of two dimensional surfaces. On a Riemannian manifold $M$, choose a pair of linearly independent vectors $u, v \in T_p M$. A two dimensional submanifold of $M$ is swept out by the set of all geodesics starting at $p$ and with initial velocities lying in the two-dimensional section $\pi$ spanned be $u, v$. The curvature of this submanifold is called the sectional curvature at $p$ of the section $\pi$.

In all subsequent sections, we assume that $M$ has all sectional curvatures bounded above by some $C \geq 0$.

The next result, due to Kendall [22], gives a sufficient condition for the existence of a unique local minimum of $F$ in a geodesic ball of reasonably wide radius.

**Proposition 4.1.** *If the support of $Q$ is contained in $B(p, r)$ with $r < \pi/(2\sqrt{C})$ and $\overline{B(p, r)} \cap C(p) = \phi$, then the Fréchet function $F$ of $Q$ has a unique local minimum in $B(p, r)$.*

Recall that (Karchar [19]; see also Theorem 2.1 in Bhattacharya and Patrangenaru [7]) if $Q(C(p)) = 0 \; \forall p \in M$, then every local minimum $\mu$ of $F$ satisfies

$$\int_{T_\mu M} v \tilde{Q}(dv) = 0 \tag{4.2}$$

where $\tilde{Q}$ is the image of $Q$ under the map $\exp_\mu^{-1}$ on $M \setminus C(\mu)$.

### 4.1. Asymptotic distribution of the sample intrinsic mean

One can use Theorem 2.3 to get the asymptotic distribution of the sample intrinsic mean. For that we need to verify assumptions (i) to (vii). The next result gives sufficient conditions for those assumptions to hold.

**Theorem 4.2.** *Suppose the support of $Q$ is contained in a geodesic ball $B(p,r)$ with center $p$ and radius $r$ as in Proposition 4.1. Let $\phi = \exp_p^{-1} : B(p,r) \longrightarrow T_pM(\approx \Re^d)$. Define $h(x,y) = d_g^2(\phi^{-1}x, \phi^{-1}y)$; $x, y \in \Re^d$. Let $((D_rh))_{r=1}^d$ and $((D_rD_sh))_{r,s=1}^d$ be the matrices of first and second order derivatives of $y \mapsto h(x,y)$. Let $\tilde{X}_j = \phi(X_j) (j=1,\ldots,n)$, $X_1,\ldots,X_n$ being i.i.d. observations from $Q$. Let $\mu = \phi(\mu_I)$, $\mu_I$ being the point of local minimum of $F$ in $B(p,r)$. Let $\mu_n = \phi(\mu_{nI})$, $\mu_{nI}$ being the point of local minimum of $F_n$ in $B(p,r)$. Define $\Lambda = E((D_rD_sh(\tilde{X}_1, \mu)))_{r,s=1}^d$, $\Sigma = \text{Cov}((D_rh(\tilde{X}_1, \mu)))_{r=1}^d$. If $\Lambda$ and $\Sigma$ are nonsingular, then*

$$(4.3) \qquad \sqrt{n}(\mu_n - \mu) \xrightarrow{\mathcal{L}} N(0, \Lambda^{-1}\Sigma\Lambda^{-1}).$$

*Proof.* When $Q$ is considered as a probability measure on the compact ball $\overline{B(p,r)}$ (as the underlying metric space), $\mu_n$ is a consistent estimator of $\mu$, by Proposition 2.2. In view of Proposition 4.1, as in the proof of Theorem 2.3 in Bhattacharya and Patrengenaru [8], Assumptions (i)–(vii) of Theorem 2.3 are verified. □

**Remark 4.1.** The nonsingularity of $\Sigma$ in Theorem 4.2 is a mild condition which holds in particular if $Q$ has a density (component) with respect to the volume measure. The nonsingularity of $\Lambda$ is a more delicate matter in general, involving a detailed analysis involving curvature and Jacobi fields. These matters are considered in detail in Bhattacharya and Bhattacharya [3].

**Remark 4.2.** Under the hypothesis of Proposition 4.1 (and Theorem 4.2), the point of local minimum $\mu_I$ of $F$ in $B(p,r)$ may not be the global minimizer of $F$ on $M$. However, if one restricts attention to the closed ball $\overline{B(p,r)}$ as the underlying metric space of interest, this point of local minimum is the intrinsic mean (on $\overline{B(p,r)}$). The advantage of Theorem 4.2 over the earlier result Theorem 2.3 in Bhattacharya and Patrengenaru [8] is that here one allows a much wider support of $Q$, namely, the radius $r$ here is twice as large as that allowed in the earlier result. This is particularly important in two-sample problems as well as in problems of classification involving several populations. Also from a statistical point of view, the mean shape is perhaps better represented if defined as the Fréchet mean over $\overline{B(p,r)}$ than over the whole of $M$, since $Q(M \setminus \overline{B(p,r)}) = 0$ and since $B(p,r)$ is a connected Riemannian manifold inheriting the metric of $M$.

Theorem 4.2 can be used to construct an asymptotic $1 - \alpha$ confidence set for $\mu_I$ which is given by

$$(4.4) \qquad \{\mu_I : n(\mu_n - \mu)^t(\hat{\Lambda}^{-1}\hat{\Sigma}\hat{\Lambda}^{-1})^{-1}(\mu_n - \mu) \leq \mathcal{X}_d^2(1-\alpha)\}$$

where $(\hat{\Sigma}, \hat{\Lambda})$ are consistent sample estimates of $(\Sigma, \Lambda)$ and $\mathcal{X}_d^2(1-\alpha)$ is the upper $(1-\alpha)^{\text{th}}$ quantile of the chi-squared distribution with $d$ degrees of freedom.

Also we can perform a nonparametric test to test if two distributions $Q_1$ and $Q_2$ have the same intrinsic mean $\mu_I$. Let $\mu = \phi(\mu_I)$. Let $X_1,\ldots,X_n$ and $Y_1,\ldots,Y_m$ be i.i.d. observations from $Q_1$ and $Q_2$, respectively. Let $Q_n$ and $Q_m$ be the empirical distributions and $\mu_{n1}$ and $\mu_{m2}$ be the corresponding sample mean coordinates. We want to test $H_0 : \mu_{1I} = \mu_{2I} = \mu_I$, say, against $H_1 : \mu_{1I} \neq \mu_{2I}$, where $\mu_{1I}$ and $\mu_{2I}$ are the true intrinsic means of $Q_1$ and $Q_2$, respectively. Then the test statistic used is

$$(4.5) \qquad T_{nm} = (n+m)(\mu_{n1} - \mu_{m2})' \hat{\Sigma}^{-1} (\mu_{n1} - \mu_{m2}),$$

$$(4.6) \qquad \hat{\Sigma} = (m+n) \left( \frac{1}{n} \hat{\Lambda}_1^{-1} \hat{\Sigma}_1 \hat{\Lambda}_1^{-1} + \frac{1}{m} \hat{\Lambda}_2^{-1} \hat{\Sigma}_2 \hat{\Lambda}_2^{-1} \right),$$

$(\Lambda_1, \Sigma_1)$ and $(\Lambda_2, \Sigma_2)$ being the parameters in the asymptotic distribution of $\sqrt{n}(\mu_{n1} - \mu)$ and $\sqrt{m}(\mu_{m2} - \mu)$, respectively, as defined in Theorem 4.2. $(\hat{\Lambda}_1, \hat{\Sigma}_1)$ and $(\hat{\Lambda}_2, \hat{\Sigma}_2)$ are consistent sample estimates. In case $n, m \to \infty$ such that $n/(m+n) \to \theta$, $0 < \theta < 1$, then under the hypothesis of Theorem 4.2, assuming $H_0$ to be true,

$$(4.8) \qquad \sqrt{n+m}(\mu_{n1} - \mu_{m2}) \xrightarrow{\mathcal{L}} N_d(0, \frac{1}{\theta} \Lambda_1^{-1} \Sigma_1 \Lambda_1^{-1} + \frac{1}{1-\theta} \Lambda_2^{-1} \Sigma_2 \Lambda_2^{-1}).$$

So $T_{nm} \xrightarrow{\mathcal{L}} \mathcal{X}_d^2$. We reject $H_0$ at asymptotic level $1 - \alpha$ if $T_{nm} > \mathcal{X}_d^2(1-\alpha)$.

We conclude with the test for the equality of intrinsic variations $V_1$, $V_2$ of $Q_1$ and $Q_2$. Under the hypothesis of Theorem 4.2, the test for $H_0 : V_1 = V_2$, against $H_1 : V_1 \neq V_2$, is provided by the asymptotically Normal statistic $T_{nm}$ in (2.14), as described at the end of Section 2.

## 5. Examples

In this section, we record the results of two-sample tests in two examples.

**Example 1** (Schizophrenic Children). In this example from Bookstein [9], 13 landmarks are recorded on a midsagittal two-dimensional slice from a Magnetic Resonance brain scan of each of 14 schizophrenic children and 14 normal children. Figures 1(a), (b) show the preshapes of the landmarks for the patient and normal samples along with the respective sample extrinsic mean preshapes. The sample preshapes are rotated appropriately as to minimize their Euclidean distance from the mean preshape. Figure 2 shows the preshapes of the normal and the patient sample extrinsic means along with the pooled sample mean.

The values of the two-sample test statistics (3.17), (4.5) for testing equality of the mean shapes, along with the p-values are as follows.

      Extrinsic: $T_{nm} = 95.5476$, p-value $= P(\mathcal{X}_{22}^2 > 95.5476) = 3.8 \times 10^{-11}$.
      Intrinsic: $T_{nm} = 95.4587$, p-value $= P(\mathcal{X}_{22}^2 > 95.4587) = 3.97 \times 10^{-11}$.

The extrinsic sample variations for patient and normal samples are 0.0107 and 0.0093, respectively. The value of the two-sample test statistic (3.19) for testing equality of extrinsic variations is 0.9461, and the p-value is 0.3441. The value of the likelihood ratio test statistic, using the so-called *offset normal shape distribution* (Dryden and Mardia [12], pp. 145–146) is $-2 \log \Lambda = 43.124$, p-value $= P(\mathcal{X}_{22}^2 > 43.124) = 0.005$. The corresponding values of Goodall's F-statistic and Bookstein's Monte Carlo test (Dryden and Mardia [12], pp. 145–146) are $F_{22,572} = 1.89$, p-value $= P(F_{22,572} > 1.89) = 0.01$. The p-value for Bookstein's test $= 0.04$.

**Example 2** (Gorilla Skulls). To test the difference in the shapes of skulls of male and female gorillas, eight landmarks are chosen on the midline plane of the skulls of 29 male and 30 female gorillas. We use the data of O'Higgins and Dryden reproduced in Dryden and Mardia ([12], pp. 317–318). The statistics (3.17) and (4.5) yield the following values:

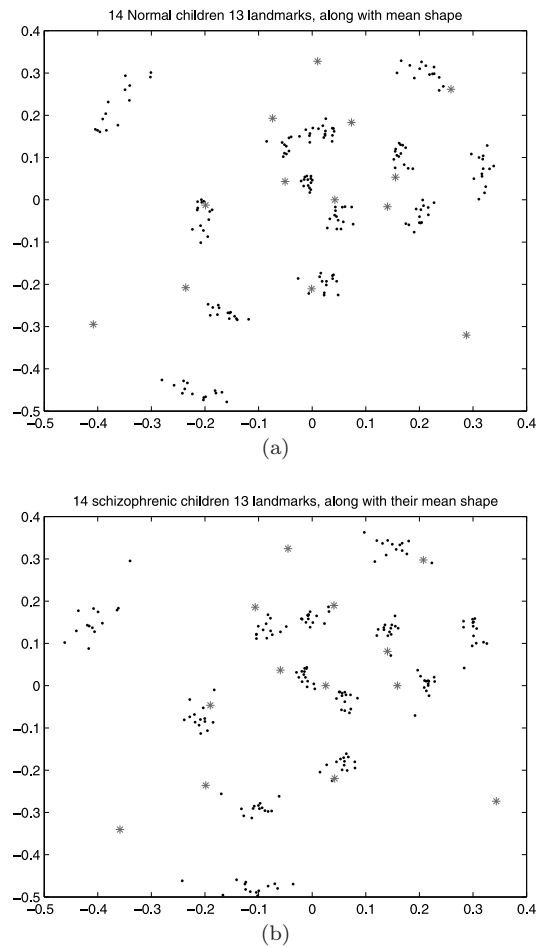

Fig 1. (a) *and* (b) *show 13 landmarks for 14 normal and 14 schizophrenic patients, respectively, along with the mean shapes, * correspond to the mean landmarks; 1c shows the sample extrinsic means for the 2 groups along with the pooled sample mean.*

$$\text{Extrinsic: } T_{nm} = 392.6, \text{ p-value} = P(\mathcal{X}_{12}^2 > 392.6) < 10^{-16}.$$
$$\text{Intrinsic: } T_{nm} = 391.63, \text{ p-value} = P(\mathcal{X}_{12}^2 > 391.63) < 10^{-16}.$$

The extrinsic sample variations for male and female samples are 0.005 and 0.0038, respectively. The value of the two-sample test statistic (3.19) for testing equality of extrinsic variations is 0.923, and the p-value is 0.356. A parametric F-test (Dryden and Mardia [12], p. 154) yields $F = 26.47$, p-value $= P(F_{12,46} > 26.47) = 0.0001$. A parametric (Normal) model for Bookstein coordinates leads to the Hotelling's $T^2$ test (Dryden and Mardia [12], pp. 170–172) yields the p-value 0.0001.

**Acknowledgments.** The authors greatly appreciate the kind and helpful suggestions by the editors and an anonymous referee which led to a substantial improvement in exposition.

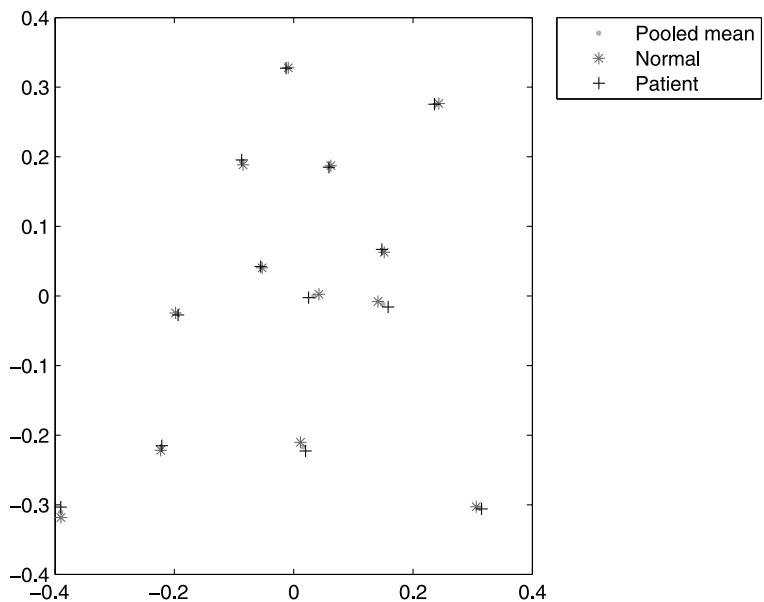

Fig 2. *The sample extrinsic means for the 2 groups along with the pooled sample mean, corresponding to Figure* 1.